\documentclass[12pt]{article}
\usepackage{amssymb,amsmath}
\def\R{\mathbb R}
\def\C{\mathbb C}
\def\Z{\mathbb Z}

\def\P{\mathbb P}
\title{Higher order Poincar\'e-Pontryagin functions and iterated path integrals}
 \author{Lubomir Gavrilov \\
 \normalsize \it Laboratoire Emile Picard, CNRS UMR 5580,  Universit\'e
 Paul Sabatier\\
 \normalsize \it 118, route de Narbonne, 31062 Toulouse Cedex, France  }
\date{November 2, 2004}
\begin{document}
\maketitle
\newtheorem{definition}{Definition}
\newtheorem{remark}{Remark}
\newtheorem{theorem}{Theorem} 
\newtheorem{lemma}{Lemma}
\newtheorem{proposition}{Proposition}
\newtheorem{corollary}{Corollary}
\vspace{5mm}
\noindent 2000 MSC scheme numbers: 34C07, 34C08
\begin{abstract}
We prove that the higher order Poincar\'e-Pontryagin functions
associated to the perturbed polynomial foliation defined by
$$df~-~\varepsilon~(P dx + Q dy )~=~0$$ satisfy a differential equation of Fuchs
type.
\end{abstract}
\section{Statement of the result} Let $f,P,Q \in \R[x,y]$ be real polynomials in
two variables. How
many limit cycles the perturbed foliation
\begin{equation} df - \varepsilon (P dx + Q dy ) = 0
\label{eq}
\end{equation} can have ? This problem is usually referred to as the weakened
16th Hilbert problem (see
Hilbert
\cite{hilbert}, Arnold
\cite[p.313]{arnold88}).

 Suppose that the foliation defined on the real plane by $\{d f = 0 \}$
possesses a family
 of periodic orbits
$\gamma (t)\subset f^{-1}(t)$, continuously depending on a parameter $t\in
(a,b)\subset \R $.
Take a segment
$\sigma
$, transversal to each orbit
$\gamma (t)$ and suppose that it can be
parameterized by $t= f|_\sigma $ (this identifies $\sigma$ to $(a,b)$).
The first return map ${\cal
P}(t,\varepsilon) $
associated to the period annulus ${\cal
A} = \cup _{t\in
\sigma } \gamma (t) $ and to (\ref{eq}) is analytic in $t,\varepsilon $ and can
be expressed
as
\begin{equation} {\cal P}(t,\varepsilon )= t + \varepsilon ^k M_k(t)+
\varepsilon ^{k+1} M_{k+1}(t)+
\dots
\label{returnmap}
\end{equation} where $M_k(t) \not\equiv 0$ is the kth order
Poincar\'e-Pontryagin function. The maximal
number of the zeros of $ M_k$ on
$\sigma $ provides an upper bound for the number of the limit cycles bifurcating
from the annulus ${\cal
A}$. For this reason $M_k(t)$ was called in \cite{gavil} {\em generating
function of limit
cycles}. The above construction can be carried out in the complex domain. In this
case the polynomials
$f,P,Q$ are complex and $\gamma (t)$ is a continuous family of closed loops
contained in the fibers
$f^{-1}(t)$, parameterized by a transversal open disc $\sigma ^{\C}$. The maximal
number of the complex
zeros of the generating function
$M_k$ on $\sigma ^{\C}$ provides an upper bound for the complex limit cycles
bifurcating from the family
$\{\gamma (t)\}_t$, see \cite{ily02}.

The main
result of the paper is the following
\begin{theorem} The generating function of limit cycles $M_k$ satisfies a linear
differential equation
of Fuchs type.
\label{main}
\end{theorem}
We show also that the monodromy group of $M_k$ is contained in $SL(n,\Z)$ where
$n$ is the order of the
equation. In this sense the differential equation satisfied by $M_k$ is of
``Picard-Fuchs" type too.
As a by-product we prove that $n \leq r^k$ where $r~=\dim H_1(f^{-1}(t_0), \Z)$
and $t_0$ is a typical value~of $f$. It is not clear, however, whether there
exists an uniform bound in
$k$ for the order $n$. In the explicit examples known to the author $n\leq r$.

 In the case $k~=~1$ the generating function $M_k$ is an Abelian integral
depending on  a
parameter
\begin{equation}
 M_1(t)= \int_{\gamma (t)} Pdx+Qdy
\label{m1}
\end{equation}
 and hence it satisfies a Fuchs equation of order at most $r$ (this bound is
exact). The identity
(\ref{m1}) goes back at least to  Pontryagin \cite{pontryagin} and has been
probably known to
Poincar\'e.
 In the case $k>1$ the (higher order) Poincar\'e-Pontryagin function $M_k$
is not necessarily of the form (\ref{m1}) with $P,Q$ rational
functions. This fact is discussed in Appendix B. We show in
section \ref{sintegral} that  $M_k(t)$ is a linear combination of
iterated path integrals of length $k$ along $\gamma (t)$ whose
entries are essentially rational one-forms.  This observation is
crucial for the proof of Theorem \ref{main}. It implies that the
monodromy representation of $M_k$ is finite-dimensional, as well
that $M_k$ is a function of moderate growth. The universal
monodromy representation of {\it all} generating functions $M_k$
was recently described in \cite{gavil}. It is not known, however,
whether this representation is finite-dimensional.

We note that iterated path integrals appeared recently in a
similar context in the study of  the polynomial Abel equation
\cite{bri03,bru03,francoise02}. Some of their basic properties
used in the paper are summarized in the Appendix.

The author acknowledges the  stimulating discussions and comments
of I.D. Iliev,  Yu. S. Ilyashenko, S. Yakovenko, and Y. Yomdin.
\section{The integral representation of $M_k(t)$.}
\label{sintegral}
From now on we consider (\ref{eq}) as a perturbed  complex foliation in $\C^2$.
Let $l (t) \in f^{-1}(t)
\subset \C^2$ be a continuous family of   closed loops, defined for all $t$
which belong to some complex
neighborhood of the typical value $t_0$ of
$f$. There exists a constant
$c>0$ such that
  the holonomy (or monodromy) map ${\cal} P(t,\varepsilon )$ of the foliation
((\ref{eq}) associated to
the family  $l (t)$ is well defined and analytic in $\{(t,\varepsilon ):
|t-t_0|<c,|\varepsilon |<c\}$.
Therefore it has there the representation (\ref{returnmap}).  Of course the
continuous deformation of a
given closed loop $l(t_0)\subset f^{-1}(t_0)$ is not unique. The free homotopy
class of the loop is
however unique and the first non-zero Poincar\'e-Pontryagin function $M_k(t)$,
 defined by (\ref{returnmap}) depends only on the free homotopy class $\gamma
(t)$ of $l(t)$ \cite{gavil}.
The main result of this section is the following
\begin{theorem}
\label{integral}
 Let $\gamma (t): [0,1] \rightarrow  f^{-1}(t)$ be a continuous family of closed
loops. For every regular
value
$t_0$ of $f$ there exists a neighborhood $U_0$ of $t_0$ in which  the first
non-vanishing
Poincar\'e-Pontryagin function $M_k(t)$,  associated to
$\gamma (t)$ and (\ref{eq}) is a finite linear combination of iterated integrals
of length at most $k$,
whose entries are differential one-forms analytic in
$f^{-1}(U_0)$.
\end{theorem}

The function $M_k(t)$ is computed according to the
 Fran\c coise's recursion formula \cite{francoise96}
$$ M_k(t) = \int_{\gamma (t)}\Omega _k
$$  where
\begin{equation}
\Omega _1= Pdx +Qdy, \Omega _m = r_{m-1} (Pdx+Qdy), 2\leq m\leq k
\label{recursion}
\end{equation}
 and the functions $r_i$ are determined successively from the (non-unique)
representation $\Omega _i = d
R_i + r_i df$. We intend to derive   explicit expressions for the functions
$r_i$. For this purpose
consider a trivial smooth
fibration
$$ f: V_0 \rightarrow U_0 = \{t\in C : |t-t_0|<c \}
$$
where $V_0$ is a connected analytic two-dimensional manifold, and $f$ is an
analytic surjection.
The fibers $f^{-1}(t)$ are mutually diffeomorphic Riemann surfaces.
Suppose that
 there exists an analytic curve
\begin{equation}
\tau : t \rightarrow P_0(t) \in f^{-1}(t)\subset V_0, t\in U_0
\label{curve}
\end{equation}
 transversal to the fibers $f^{-1}(t)$. For an analytic one-form in $V_0$ define
the function
$$F(P)=\int_{P_0(t)}^P \omega $$
 where
$t= f(P)$ and the integration is along some path contained in $f^{-1}(t)$ and
connecting the points
$P_0(t), P \in f^{-1}(t)$. Finally we shall suppose that when varying $P\in
f^{-1}(U_0)$ the path
connecting
$P_0(t)$ and $P$ varies continuously in $P$. The function $F(P) $  is
multivalued but locally  analytic
in
$V_0=f^{-1}(U_0)$.
\begin{lemma}
\label{idl}
Under the above conditions the following identity holds
\begin{equation}
\label{ide}
d \int_{P_0(t)}^P \omega = (\int_{P_0(t)}^P \frac{ d\omega}{df}) df + \omega  -
(\tau
\circ f)^* \omega
\end{equation}
where $\frac{ d\omega}{df}$ is the Gelfand-Lerray form of $d \omega $ and $
(\tau \circ
f)^* \omega$ is the pull back of $\omega $ under the map
$$
\tau \circ f : V_0 \stackrel{f}{\rightarrow} U_0 \stackrel{\tau }{\rightarrow}
V_0  \; .
$$
\end{lemma}
{\bf Remark.} If $\tilde{\tau } : t\rightarrow  \tilde{P}_0(t) \in f^{-1}(t)$ is
another transversal
curve (as in (\ref{curve})) then (\ref{ide}) implies
$$
d \int_{P_0(t)}^{\tilde{P}_0(t)}\omega = (\int_{P_0(t)}^{\tilde{P}_0(t)} \frac{
d\omega}{df}) df +
(\tilde{\tau}
\circ f)^* \omega  - (\tau
\circ f)^* \omega .
$$
If, in  particular
 $P\equiv P_0(t)$ (so the path of integration $\tau $ is closed) we get the well
known
identity \cite{avg}
$$
d \int_{\tau  (t)} \omega = (\int_{\tau  (t)} \frac{ d\omega}{df}) dt .
$$
 {\bf Proof of Lemma \ref{idl}.} Suppose that when $t$ is sufficiently close to
$t_0$ the path of integration connecting
$P_0(t)$ to $P$ is contained in some open poly-disc
$D _\varepsilon \subset V_0$, in which we may choose local coordinates $x,y$.
We claim that for such a
family of paths and for $t$ sufficiently close to $t_0$ the identity
$(\ref{ide})$ holds true. As
(\ref{ide}) is linear in $\omega $ we may suppose without loss of generality
that
 $\omega|_{D_\varepsilon} = Q(x,y) dx$ where $Q$ is analytic in $D_\varepsilon$.
Suppose further that the
path of integration from
$P_0(t)=(x_0(t),y_0(t))$ to
$P=(x,y)  $  is projected under the map $(x,y) \rightarrow x$ into an analytic
path, avoiding the
ramification points of this projection, and connecting
$x_0(t)$ and
$x$ in the complex
$x$-plane. Along such a path we may express $y=y(x,t)$ from the identity
$f(x,y)=t$ and hence
\begin{eqnarray*}
d \int_{P_0(t)}^P \omega & = & d  \int_{x_0(t)}^x Q dx \\ & = & Q dx - Q x_0'(t)
df +
(\int_{x_0(t)}^x Q_y
\frac{\partial y}{\partial t} dx) df \\
& = & Q dx - Q x_0'(t) df + (\int_{x_0(t)}^x  \frac{Q_y}{f_y} dx)
df \\
& = &   \omega  - (\tau \circ f)^* \omega + (\int_{P_0(t)}^P \frac{
d\omega}{df}) df .
\end{eqnarray*}
Therefore (\ref{ide}) holds true in a neighborhood of $P_0(t_0)$. By analytic
continuation it holds true for arbitrary $P$ and arbitrary continuous  family of
paths  connecting
$P_0(t)$ to $P$.
 The Lemma is proved.$\triangle$\\
 Let $f\in \C[x,y]$,
$\gamma (t)\subset f^{-1}(t)$ be a continuous family of closed loops such that
$\int_{\gamma (t)} \omega
\equiv 0$. If $\gamma (t)$ generates the fundamental group of $f^{-1}(t)$  then
(\ref{ide}) implies that
$\omega = dA + Bdf$ where $A,B$ are  analytic functions in $f^{-1}(U_0)$. In the
case when the
fundamental group of $f^{-1}(t)$ is not infinite cyclic we consider a covering
\begin{equation} \tilde{V_0} \stackrel{p }{\rightarrow} V_0
\label{cov}
\end{equation} such that the fundamental group of $\tilde{V_0}$ is infinite
cyclic with a generator
represented by a closed loop
$\tilde{\gamma }(t_0)$ projected to
$\gamma (t_0)$ under
$p$. Such a covering exists and is unique up to an isomorphism \cite{ful95}.
Moreover $\tilde{V_0}$ has
a canonical structure of analytic two-manifold induced by
$p$. If we define $\tilde{f}= f\circ p$, then the fibration
\begin{equation}
\tilde{f}: \tilde{V_0} \rightarrow U_0
\label{fibration}
\end{equation} is locally trivial and the fibers are homotopy equivalent to
circles. An analytic function
(or differential form) on
$\tilde{V_0}$ is a locally analytic function (differential form) on
$V_0=f^{-1}(U_0)$ such that
\begin{description}
\item[(i)] it has an analytic continuation along any arc in $f^{-1}(U_0)$
\item[(ii)] its determination does not change as $(x,y)$ varies along any closed
loop homotopic to
$\gamma (t_0)$.
\end{description} We shall denote the space of such functions (differential
forms) by $\tilde{{\cal
O}}(f^{-1}(U_0))$  ($\tilde{\Omega ^k}(f^{-1}(U_0))$). Lemma \ref{idl} implies
the following
\begin{corollary}  If $\tilde{\omega} \in \tilde{\Omega}^1(f^{-1}(U_0))$ is such
that
$\int_{\gamma (t)} \tilde{\omega} \equiv 0$, then
$$
\tilde{\omega}  = d \tilde{A} + \tilde{B} df
$$ where $\tilde{A},\tilde{B} \in \tilde{{\cal O}}(f^{-1}(U_0))$,
$$
\tilde{A}=   \int_{P_0(t)}^P \tilde{\omega} , \tilde{B}= - \int_{P_0(t)}^P
\frac{ d\tilde{\omega}}{df} +
R(f)
$$ and $R(.)$ is analytic in $U_0$.
\label{dif}
\end{corollary}
In the proof of  Theorem \ref{integral} we shall use the following well known
\begin{proposition}
\label{grad}
 Let $f\in \C[x,y]$ be a non-constant polynomial. Then there exists a polynomial $m\in \C[f]$ such that
\begin{enumerate}
    \item $m(f)$ belongs to the gradient ideal of $f$
    \item $m(c)=0$ if and only if $c$
is a critical value of $f$.
\end{enumerate}
\end{proposition}
The identity
\begin{equation}
\label{grad1} (\alpha f_x + \beta f_y ) dx \wedge dy = df \wedge
(\alpha dy - \beta dx)
\end{equation}
combined with Proposition \ref{grad} shows that when $\omega$ is a
polynomial (analytic) one-form, then the Gelfand -Leray form
$$
m(f) \frac{d \omega }{df}
$$
can be chosen polynomial (analytic).\\
{\bf Proof of Proposition \ref{grad}} Consider the reduced
gradient ideal $J_{red}\subset \C[x,y]$ generated by $f_x/D,
f_y/D$ where $D$ is the greatest common divisor of $f_x,f_y$. The
variety $V(J_{red})= \{ c_i \}_i$ is a finite union of points
which may be supposed non-empty. Therefore $\C[x,y]/J_{red}$ is a
vector space of finite dimension\cite{ful95} and the
multiplication by $f$ defines an endomorphism. Therefore $m(f) \in
J_{red}$ where $m(.)$ is the minimal polynomial of the
endomorphism defined by $f$. We note that $m(c)=0$ if and only if
$c=c_i$ for some $i$. Taking into consideration that $\prod_i
(f-c_i)/D$ is a polynomial we conclude that $ \prod_i (f-c_i) m(f)
$ belongs to the gradient ideal of $f$.\\
{\bf Proof of Theorem \ref{integral}.}
 Suppose that $M_1=...=M_{k-1}=0$ but $M_k \neq 0$, $k\geq 3$. The recursion
formula (\ref{recursion})
implies that
$M_i(t) = \int_{\gamma (t)} \tilde{\omega}_i$, where $\tilde{\omega }_i \in
\tilde{{\Omega
}}^1(f^{-1}(U_0)) $, $i\leq k$. Indeed, $\tilde{\omega }_1=\omega \in {\Omega
}^1(f^{-1}(U_0))$ and if
$\tilde{\omega }_i \in \tilde{{\Omega }}^1(f^{-1}(U_0)) $, then by Corollary
\ref{dif}
$$\tilde{\omega }_{i+1} = -  \omega \int_{P_0(t)}^P \frac{d \tilde{\omega
}_i}{df} .
$$ The Gelfand-Leray form $\frac{d \tilde{\omega }_i}{df}$ may be supposed
analytic (according to
Proposition \ref{grad} and (\ref{grad1}) ).
$M_i'(t)= \int_{\gamma (t)}\frac{d \tilde{\omega }_i}{df}=0$ implies that
$\int_{P_0(t)}^P \frac{d \tilde{\omega }_i}{df} \in \tilde{\cal O
}(f^{-1}(U_0))$ and hence
$\tilde{\omega }_{i+1} \in \tilde{{\Omega }}^1(f^{-1}(U_0))$. We obtain in
particular that
$$ M_k(t) = - \int_{\gamma (t)} \omega \int_{P_0(t)}^P \frac{d \tilde{\omega
}}{df}
$$ where
$$ M_{k-1}(t)= \int_{\gamma (t)}  \tilde{\omega } \equiv 0 .
$$ We shall prove the Theorem by induction on $k$. Suppose that that
$M_{k-1}(t)$,is a finite linear
combination of iterated integrals of  length at most $k-1$, whose entries are
differential one-forms
analytic in
$f^{-1}(U_0)$. We need to show that the same holds true for
\begin{equation}
\int_{P_0(t)}^P \frac{d \tilde{\omega }}{df} .
\label{eqd}
\end{equation} Let $\omega _1,\omega _2,\dots ,\omega _{k-1}$ be analytic
one-forms in $f^{-1}(U_0)$.
Lemma \ref{idl} implies
$$
\frac{d }{df} (\omega _1 \int_{P_0(t)}^{P_1} \omega _2 \int_{P_0(t)}^{P_2}
\omega _3 \dots
\int_{P_0(t)}^{P_{k-1}} \omega _{k-1} ) =
\frac{d
\omega _1 }{df}  \int_{P_0(t)}^{P_1} \omega _2  \dots  \int_{P_0(t)}^{P_{k-1}}
\omega _{k-1}
$$
$$ +
\omega _1   \int_{P_0(t)}^{P_1} \frac{d \omega _2}{df}  \dots
\int_{P_0(t)}^{P_{k-1}} \omega _{k-1} +
\dots  + \omega_1
\int_{P_0(t)}^{P_1} \omega _2  \dots \int_{P_0(t)}^{P_{k-1}} \frac{d \omega
_{k-1}}{df}
$$
$$ -
\frac{\omega _1\wedge \omega _2}{df} \int_{P_0(t)}^{P_2} \omega _3 \dots
\int_{P_0(t)}^{P_{k-1}} \omega
_{k-1}
$$
$$ -\omega_1 \int_{P_0(t)}^{P_1} \frac{\omega _2\wedge \omega _3}{df}
\int_{P_0(t)}^{P_3} \omega _4 \dots
\int_{P_0(t)}^{P_{k-1}}
\omega _{k-1}
$$
$$
\dots
$$
$$ - \omega _1 \int_{P_0(t)}^{P_1} \omega _2 \dots \int_{P_0(t)}^{P_{k-4}}
\omega
_{k-3}\int_{P_0(t)}^{P_{k-3}}
\frac{\omega _{k-2} \wedge \omega _{k-1}}{df}
$$
$$ + \omega _1 \int_{P_0(t)}^{P_1} \omega _2 \int_{P_0(t)}^{P_2} \omega _3 \dots
 \int_{P_0(t)}^{P_{k-2}}
\frac{\omega _{k-2}\wedge (
\tau
\circ f)^* \omega _{k-1}}{df} .
$$ The differential form $\frac{\omega _{k-2}\wedge ( \tau
\circ f)^* \omega _{k-1}}{df}$ can be written in the form $\omega _{k-2} R(f)$
where
$(\tau
\circ f)^* \omega _{k-1} = -R(f) df$.  This shows that $(\ref{eqd})$  is a
linear combination of iterated
integrals of length at most $k-1$. As the Gelfand-Leray forms
$$
\frac{\omega _i\wedge \omega _{i+1}}{df}, \frac{d \omega i}{df}
$$ may always be chosen analytic in $f^{-1}(U_0)$ (see (\ref{grad1})) then
 Theorem \ref{integral} is proved. $\Box$\\  The proof of the above theorem
provides also an algorithm
for computing the higher-order Poincar\'e-Pontryagin functions $M_k$ in terms of
iterated integrals.  To
illustrate this we consider  few examples. To simplify the notations, for every
given one-form
$\omega $ on $\C^2$, we denote by
$\omega '$ some {\it fixed} one-form, such that $df\wedge w'= d\omega $ (that is
to say $\omega '$ is a
Gelfand-Leray form of
$d \omega $).
\\ {\bf Examples}
\begin{enumerate}
\item It is well known that
$$ M_1(t) = \int_{\gamma (t)} \omega .
$$
\item If $M_1 =0$ then  Lemma \ref{idl} and Corollary \ref{dif} imply
$$
\omega = dA(x,y) + B(x,y) df + d R(f)
$$ where
$$
 A(x,y) = \int_{P_0(t)}^P \omega  , B(x,y)= - \int_{P_0(t)}^P \frac{ d \omega
}{df} + R(f) , P=(x,y) .
$$
 We have
$$  M_2(t) = \int_{\gamma (t)} B \omega = - \int_{\gamma (t)} \omega
\int_{P_0(t)}^{P} \omega '
$$  and hence
\begin{equation}
\label{m2} M_2(t) = - \int_{\gamma (t)}  \omega
\omega '.
\end{equation} As $ \int_{\gamma (t)}  \omega ' = \int_{\gamma (t)}  \omega
\equiv 0$ then the  iterated
integral (\ref{m2}) depends on the free homotopy class of $\gamma (t)$ (and not
on the initial point
$P_0(t)$ ).
\item If $M_2=0$, then
$$ M_3(t)= \int_{\gamma (t)} \omega  \int_{P_0(t)}^{P}\frac{d}{df}( \omega
\int_{P_0(t)}^{Q} \omega ')
$$ where
$$
\frac{d}{df} (\omega \int_{P_0(t)}^{P} \omega ')=
\omega ' \int_{P_0(t)}^{P} \omega ' + \omega \int_{P_0(t)}^{P} \omega '' -
\int_{\gamma (t)}\frac{\omega
\wedge \omega '}{df} +R(t) \int_{P_0(t)}^{P} \omega
$$ and $R(t)$ is an analytic function computed from the identity $ (f\circ \tau
)^* \omega '= R(f)df$.
As $\int_{\gamma (t)} \omega
\omega =0 $ then
\begin{equation}
 M_3(t)=
\int_{\gamma (t)} \omega (\omega ')^2  + \int_{\gamma (t)}\omega ^2\omega ''-
 \omega  \frac{\omega \wedge \omega '}{df}  .
\label{m3}
\end{equation} Both of the iterated integrals in (\ref{m3}) depend on the free
homotopy class of $\gamma
(t)$ only and do not depend on the particular choice of the Gelfand-Leray form
$\omega '$.
\item If $M_3=0$, then
\begin{eqnarray*} M_4(t)& = &  - \int_{\gamma (t)} \omega
\int_{P_0(t)}^{P}\frac{d}{df} (\omega
\int_{P_0(t)}^{Q}  (\omega ')^2  )\\ & &
  - \int_{\gamma (t)} \omega  \int_{P_0(t)}^{P}\frac{d}{df} (\omega
\int_{P_0(t)}^{Q} \omega \omega ''
)\\ & &
 - \int_{\gamma (t)} \omega  \int_{P_0(t)}^{P}\frac{d}{df} (  \omega
\int_{P_0(t)}^{Q} \frac{\omega
\wedge \omega '}{df})
\end{eqnarray*} If we make the particular choice $\omega '= -d B$, $B =
\int_{P_0(t)}^{P} \frac{d \omega
}{df} $ for the Gelfand-Leray form of
$d\omega $, as well $w''=0$ then the formula for $M_4$ becomes
\begin{eqnarray*} M_4(t) & = &
\int_{\gamma (t)} \omega (\omega ')^3  + \omega \omega' \frac{\omega \wedge
\omega '}{df} +
\omega^2  \frac{d \frac{\omega \wedge \omega '}{df}}{df} \\ & & + \int_{\gamma
(t)} \omega  \frac{\omega
\wedge \omega '}{df}
\omega ' -   \omega \frac{\omega  \wedge
\frac{\omega \wedge \omega '}{df}}{df} \\ & &  - (\int_{\gamma (t)}  \omega
^2\omega ' )(f\circ \tau )^*
\omega ' + (\int_{\gamma (t)}\omega ^2)(f\circ \tau )^* \frac{\omega \wedge
\omega '}{df}\\ &=&
\int_{\gamma (t)} \omega (\omega ')^3  + \omega \omega' \frac{\omega \wedge
\omega '}{df} +
\omega^2  \frac{d \frac{\omega \wedge \omega '}{df}}{df} \\ & & + \int_{\gamma
(t)} \omega  \frac{\omega
\wedge \omega '}{df}
\omega ' -   \omega \frac{\omega  \wedge
\frac{\omega \wedge \omega '}{df}}{df} .
\end{eqnarray*} Note that the last expression depends on the choice of $\omega
'$ and hence on the
initial point $P_0(t)$. It is an open question to find a general closed formula
for $M_k$, $k\geq 4$, in
terms of iterated integrals with rational entries, depending on the free
homotopy class of $\gamma (t)$
only.
\end{enumerate}

\section{Proof of Theorem \ref{main}}
\label{s3}
The proof is split in two parts. First we show that
$M_k(t)$ satisfies a linear differential equation of finite order (possibly with
irregular
singularities). For this we need to study the monodromy group of  $M_k(t)$.
Second, we shall show that
the generating function $M_k(t)$ is of moderate growth on the projective plane
$\C\P^1$, and hence the
equation is Fuchsian.
\subsection{ The monodromy representation of $M_k$.}
\label{monrep}
 Recall first that the universal monodromy
representation for $M_k(t)$ (for arbitrary $k$) can be constructed
as follows (see \cite{gavil} for proofs).  To the non-constant
polynomial $f\in \C[x,y]$ we associate the locally trivial
fibration
$$ f^{-1}(\C \backslash \Delta ) \stackrel{f}{\rightarrow } \C \backslash \Delta
$$ where $\Delta \subset \C $ is the finite set of atypical values. Let $t_0\not \in \Delta $ and put
$S = f^{-1}(t_0)$. The canonical group homomorphism
\begin{equation}
\label{class} \pi_1(\C\setminus \Delta  , t_0) \rightarrow
\mbox{\rm Diff}\,(S )/\mbox{\rm Diff}_0(S)  .
\end{equation} where ${\rm Diff}\,(S )/\mbox{\rm Diff}_0(S ) $ is the mapping class group of $S $,
induces
 a homomorphism (group action on
$\pi_1(S)$)
\begin{equation}
\label{action1} \pi_1(\C\setminus D , t_0) \rightarrow \mbox{\it
Perm}\,(\pi_1(S ))
\end{equation} where $\mbox{\it Perm}\,(\pi_1(S ))$ is the group of permutations of
$\pi_1(S)$, and $\pi_1(S)$ is the set of free homotopy classes of
closed loops on $S$.

Let $\gamma (t_0) \in \pi_1(S) $ be a free homotopy class of
closed loops on $S$ and consider the orbit ${\cal O}_{\gamma
(t_0)}$ of $\gamma (t_0)$ under the group action (\ref{action1}).
For a given point $P_0\in S$ we denote $F= \pi_1(S,P_0)$ and let
 $G \subset \pi_1(S ,P_0)$  be the normal subgroup generated by the pre-image of the orbit
${\cal O}_{\gamma (t_0)}$ under the canonical projection
$$
\pi_1(S ,P_0)  \rightarrow \pi_1(S ) \;.
$$
 Let $(G,F)$ be the normal sub-group of $G$ generated by commutators 
$$
g^{-1}f^{-1}gf, g\in G, f\in F
$$
 and denote
 $$
 H_1^\gamma(S,\mathbb{Z})= G/(G,F) .
 $$
From the definition of $G$ it follows that the Abelian group
$H_1^\gamma(S,\mathbb{Z})$ is invariant under the action of
$\pi_1(\C\setminus \Delta )$ and hence we obtain a homomorphism
(the universal representation)
\begin{equation}
\label{monodromy1} \pi_1(\C\setminus \Delta  , t_0) \rightarrow
Aut(H_1^\gamma(S,\mathbb{Z})) .
\end{equation}
On the other hand, the monodromy representation of the generating
function
 $M_k(t) = M_k(\gamma,{\cal F}_\varepsilon,t)$ is defined as follows. The
function $M_k(t)$ is multivalued on  $\C~\setminus~\Delta $. Let
us  consider all its possible determinations in a sufficiently
small neighborhood of $t=t_0$. All integer linear combinations of
such functions form a module over $\Z$ which we denote by ${\cal
M}_k( \gamma,{\cal F}_\varepsilon )$. The fundamental group $\pi
_1(\C \setminus \Delta ,t_0) $ acts on $ {\cal M}_k ={\cal M}_k(
\gamma,{\cal F}_\varepsilon )$ in an obvious way. We obtain thus a
homomorphism
\begin{equation}
\label{repr} \pi _1(\C \setminus \Delta ,t_0) \rightarrow
Aut({\cal M}_k)
\end{equation}
called the monodromy representation of the generating function
$M_k(\gamma,{\cal F}_\varepsilon,t)$.

It is proved in \cite[Theorem 1]{gavil} that the map
\begin{equation}\label{hom}
H^{\gamma}_1(f^{-1}(t_0) ,\Z) \stackrel{\varphi}{\rightarrow}
{\cal M}_k( \gamma,{\cal F}_\varepsilon ): \gamma \rightarrow
M_k(\gamma,{\cal F}_\varepsilon,t)
\end{equation}
is a canonical surjective homomorphism compatible with the action
of the fundamental group $\pi_1(\C\setminus \Delta ,t_0)$.
Equivalently, (\ref{repr}) is a sub-representation of the
representation dual to the universal representation
(\ref{monodromy1}).

If the rank of the Abelian group $H^{\gamma}_1(f^{-1}(t_0) ,\Z)$
were finitely generated, then this would imply that each $M_k(t)$
satisfies a linear differential equation whose order is bounded by
the dimension of $H^{\gamma}_1(f^{-1}(t_0) ,\Z)$.
 We shall prove here a weaker statement:
$M_k(t)$ satisfies a linear equation of finite order (depending on
$k$). Our argument is based on the integral representation for
$M_k(t)$ obtained in the previous section. Namely, consider the
lower central series
\begin{equation}
 {F_1} \supseteq {F}_2  \supseteq \dots {F_k} \supseteq \dots
\label{central}
\end{equation}
 where $F_1= F= \pi_1(S,P_0)$, and $F_{k+1}= (F_k,F)$ is the subgroup of $F_k$ generated by commutators
$(f_k,f)=f_k^{-1}f^{-1}f_k f$, $f_k \in F_k$, $f\in F$. An
iterated integral of length $k$ along a closed loop which belongs
to $F_{k+1}$ vanishes. Therefore Theorem \ref{integral}  implies
that in (\ref{monodromy1}) we can further truncate by $F_{k+1}$.
Namely, for every subgroup $H\subset F$ we denote $\tilde{H}= (H
\cup F_{k+1}) / F_{k+1}$.
 As before the group action $(\ref{action1})$  induces a homomorphism
\begin{equation}
\label{monodromy2} \pi_1(\C\setminus D , t_0) \rightarrow
Aut(\tilde{G}/(\tilde{G},\tilde{F}))
\end{equation}
and there is a canonical surjective homomorphism
$$
\tilde{G}/(\tilde{G},\tilde{F}) \rightarrow {\cal M}_k(
\gamma,{\cal F}_\varepsilon ) .
$$
The lower central series of $\tilde{F}= \tilde{F_1}$ is
$$ \tilde{F_1} \supseteq \tilde{F}_2  \supseteq \dots \tilde{F_k} \supseteq \{id\} .
$$ It is easy to see that in this case $(\tilde{G},\tilde{F}) \subset \tilde{F}
$ is finitely generated (e.g. \cite[ Lemma 4.2, p.93]{passman}),
and hence $\tilde{G}/(\tilde{G},\tilde{F})$ is finitely generated
too. This implies on its hand that ${\cal M}_k( \gamma,{\cal
F}_\varepsilon )$ is finite-dimensional, and hence the generating
function  satisfies a linear differential equation. Its order is
bounded by the dimension  of $\tilde{G}/(\tilde{G},\tilde{F})$.
The latter is easily estimated to be less or equal to
$$
\sum_{i=1}^k \dim F_i/F_{i+1} \leq r^k
$$ (for the last inequality see \cite[section 11]{hall}).
To resume, we proved that {\it the generating function of limit
cycles $M_k(t)$, $t\in \C\setminus \Delta $, satisfies an analytic
linear differential equation of order at most $r^k$}.

\subsection{The moderate growth of $M_k$.}
We shall show that the
possible singular points (contained in $\Delta
\cup\infty$) are of Fuchs type. A necessary and sufficient condition for this is the moderate growth of
$M_k(t)$ in any sector centered at a singular point. For this we shall use once again the integral
representation for $M_k(t)$. Let $t_0$ be an atypical value for $f$ and suppose that the analytic curve
$$
\tau : t \rightarrow P_0(t) \in f^{-1}(t)
$$
 is defined for $t\sim t_0$ and is transversal to the fibers $f^{-1}(t)$. It follows from the proof of
Theorem \ref{integral}  that
$M_k(t)$ has an integral representation in a punctured neighborhood of $t_0$ too. More precisely $M_k(t)$
is a finite linear combination
$$
\sum \frac{\alpha _i(t)}{m_i(t)} \int_{\gamma (t)} \omega _1^{i}\dots \omega _{i_j}^i
$$ where
$\omega ^p_q$ are polynomial one-forms, $i_j \leq k$, $\alpha_i(t)$ are analytic functions and $m_i(t)$
are polynomials.  Let $$t \rightarrow P_0(t)\in f^{-1}(t), t \rightarrow  \tilde{P}_0(t)\in f^{-1}(t)$$
be two analytic curves defined in a neighborhood of the atypical value $t_0$ and transversal to the
fibers $f^{-1}(t)$  (including
$f^{-1}(t_0)$) and consider the iterated integral
$$F(t) = \int_{l(t)} \omega _1 \omega _2 \dots  \omega _k$$ where $l(t)$ is a path on $f^{-1}(t)$
connecting $P_0(t)$ and $\tilde{P}_0(t)$, and
 $\omega _1,\omega _2,\dots,\omega _k$ are polynomial one-forms in $\C^2$.
\begin{proposition}
\label{growth} Let $S(t_0)= \{ t\in \C: arg (t-t_0) < \varphi _0, 0<|t-t_0|< r_0 \}$ be a sector centered
at $t_0$. There exists $r_0 , N_0>0$ such that
$|F(t)|< |t-t_0|^{-N_0}$. Let $S(\infty)= \{ t\in \C: arg (t) < \varphi _0, |t|> r_0 \}$ be a sector
centered at
$\infty$. There exist $r_0, N_0 > 0$ such that  $|F(t)| < |t|^{N_0}$ in $S(\infty)$.
\end{proposition} {\bf Remark.} Recall that an analytic function $F$ defined on the universal covering of
$\C\setminus \Delta $ and satisfying the claim of the above Proposition is said to be of {\it moderate
growth}.\\
 {\bf Proof of Proposition \ref{growth}.} Let $(x_i(t),y_i(t))\in f^{-1}(t)$ be the ramification points of
the projection
$f^{-1}(t) \rightarrow \C$ induced by
 $\pi :(x,y)\rightarrow x$. Each ramification point $x_i(t)$ has a Puiseux expansion in a neighborhood of
$t_0$. Therefore when $t$ tends to $t_0$ in a sector centered at $t_0$, each ramification point tends to
a definite point $P \in \C\P^1$.

Assume further that the projection of $l(t)$ on the $x$-plane  is represented by a piece-wise straight
line
\begin{equation}
\label{pi}
\pi (l(t)) = \cup_{i=1}^n [x_i,x_{i+1}]
\end{equation} connecting $x_0(t),x_1(t),...,x_{n+1}(t)$ where $x_i$, $i= 1,2,\dots ,n$ are some
ramification points, and
$P_0(t)=(x_0(t),y_0(t))$, $\tilde{P}_0(t)=(x_{n+1}(t),y_{n+1}(t))$. The iterated integral  $F(t)$ along
$l(t)$ is expressed as an iterated integral along $\pi (l(t))$ whose entries are one-forms with algebraic
coefficients. It is clear that such an iterated integral along $[x_i(t),x_{i+1}(t)]$ is of moderate
growth (because $x_i(t)$ are of moderate growth). Thus, if the number $n$ in (\ref{pi}) were bounded
when
$t$ varies in $S(t_0)$, then $F(t)$ would be of moderate growth of $F(t)$. It remains to
show that the number $n = n(t)$ is uniformly bounded in $S(t_0)$. Note that when $t\in S(t_0)$ and $r_0$
is sufficiently small the ramification points $x_i(t)$ are all distinct. Denote by {\bf B} the subset of
$ S(t_0)$ of points $t$ having the property

``there exist ramification points $ x_i(t), x_j(t), x_k(t)$ such that $x_i(t)- x_j(t)$ and $x_i(t)-
x_k(t)$ are collinear but not identically collinear ."

Using the fact that $x_i(t)- x_j(t)$,  $x_i(t)- x_k(t)$  have Puiseux expansions, we conclude that  {\bf
B} is a real analytic subset of
$\R^2 \simeq \C$ of co-dimension one.  The set
$S(t_0) \setminus {\bf B}$ has a finite number of connected components and on each connected component
the function $t\rightarrow n(t)$ is constant. It follows that $n(t)$ is uniformly bounded.
 Proposition \ref{growth}, and hence Theorem \ref{main} is proved.

\appendix
\section{ Iterated Path Integrals}
 Let $S$ be a Riemann surface and $\omega _1, \omega _2,\dots, \omega _k$ be
holomorphic one-forms. For
every smooth path
$l : [0,1] \rightarrow S$ we define the {\em iterated path integral}
\begin{equation}
\label{int2}
\int_l \omega _1 \omega _2 \dots  \omega _k = \int_{0\leq t_k\leq \dots \leq t_1 \leq 1} f_k(t_k)
\dots  f_1(t_1) dt_k \dots dt_1
\end{equation} where
 $l^*\omega _i = f_i(t) dt$. We have for instance
$$
\int_l \omega _1 \omega _2 = \int_l \omega _1 \int_{l(0)}^{l(t)} \omega _2 .
$$
The basic properties of the iterated path integrals (\ref{int2}) were established  by Parsin \cite{par69}.
The general theory of iterated integrals has been developed by Chen, e.g.
\cite{chen1,chen2}. In the Chen's theory the iterated integrals generate the De Rham complex of the path
space $PX$ associated to an arbitrary manifold
$X$. In this context the iterated integrals of the form  (\ref{int2}) provide the 0-cochains of the path
space $P_{a,b}S$, where $a,b$ are the two ends of $l$. Indeed, a connected component of $P_{a,b}S$
consists of those paths (with fixed ends) which are homotopy equivalent, and (\ref{int2}) is constant
on such paths.
Some basic properties of iterated path integrals are summarized below. The
missing proofs (and much
more) may be found in  R. Hain \cite{hain87,hain02}.
\begin{lemma} \mbox{ }
\begin{itemize}
\item[(i)] The value of $\int_l \omega _1 \omega _2 \dots  \omega _k $
 depends only on the homotopy class of $l$ in the set of loops with ends fixed at $l(0), l(1)$.
\item[(ii)] If $l_1,l_2: [0,1] \rightarrow S$ are composable  paths (i.e. $l_1(1)=l_2(0)$ ) then
\begin{equation}
\int_{l_1l_2}  \omega _1 \omega _2 \dots  \omega _k = \sum_{i=0}^k
\int_{l_2} \omega _1 \omega _2 \dots \omega _i \int_{l_1} \omega
_{i+1} \omega _2 \dots  \omega _k \label{2ii}
\end{equation}
 where we set $ \int_{l} \omega _1 \omega _2 \dots  \omega _i = 1$ if $i=0$.
\item[(iii)]
$$
\int_l \omega _1 \omega _2 \dots  \omega _k = (-1)^k \int_{l^{-1}} \omega _k \omega _{k-1} \dots  \omega
_1 .
$$
\end{itemize}
\label{l1}
\end{lemma}
From now on we suppose that $l(0)=l(1)=P_0$ and put $F=\pi _1(S,P_0)$.
For $\alpha ,\beta \in F$ we denote the commutator $\alpha ^{-1}\beta ^{-1}\alpha \beta$ by $(\alpha
,\beta )$. If $A,B\subset F$ are subgroups, we denote by $(A,B)$ the subgroup of $F$ generated by all
commutators $(\alpha ,\beta )$, such that $\alpha \in A$ and
$\beta \in B$.
 Consider the lower central series
$ F=F_1 \supseteq F_2 \supseteq F_3  \supseteq \dots
$ where $ F_k= (F_{k-1},F)$ and $F_1=F$.
\begin{lemma} \mbox{ }
\begin{itemize}
\item[(i)]
$$
- \int_{(\alpha, \beta) } \omega _1 \omega _2 = \det \left(
\begin{array}{lr}
\int_\alpha  \omega _1 & \int_\beta \omega _1\\
\int_\alpha \omega _2 & \int_\beta \omega _2
\end{array}
\right), \forall \alpha ,\beta  \in F_1 .
$$
\item[(ii)]
 Let $\gamma \in F_k$ and $\omega _1, \omega _2,\dots,\omega _k$ be holomorphic one-forms on $S$. Then
the iterated path integral
$\int_{\gamma }~\omega _1~\omega _2~\dots~\omega _{k-1}$ vanishes, and the value of the integral
$\int_{\gamma }~\omega _1~\omega _2~\dots~\omega _{k}$   does not depend on the initial point $P_0$.
\item[(iii)] If $\alpha ,\beta \in F_k$ then
$$
\int_{\alpha \beta }~\omega _1~\omega _2~\dots~\omega _{k} =\int_\alpha ~\omega _1~\omega _2~\dots~\omega
_{k}+\int_\beta ~\omega _1~\omega _2~\dots~\omega _{k} .
$$
\item[(iv)] If $\alpha \in F_p$, $\beta \in F_q$, then
\begin{eqnarray*} - \int_{(\alpha, \beta) }~\omega _1~\omega _2~\dots~\omega _{p+q}&= &\int_{\alpha
}~\omega _1~\omega _2~\dots~\omega _{p}\int_{\beta  }~\omega _{p+1}~\omega _{p+2}~\dots~\omega _{p+q}\\
&- & \int_{\beta  }~\omega _1~\omega _2~\dots~\omega _{q}\int_{\alpha }~\omega _{q+1}~\omega
_{q+2}~\dots~\omega _{p+q}
\end{eqnarray*}
\end{itemize}
\label{l2}
\end{lemma}
{\bf Proof.}
The identity $(i)$ follows from (\ref{2ii}). It implies in particular that $\int_\gamma \omega _1\omega
_2$ does not depend on the initial point $P_0$ provided that $\gamma \in F_2$, and vanishes provided
that $\gamma \in F_3$.
Suppose that the claim (ii) is proved up to order $k-1$ and let $\gamma \in F_k$. If
$\gamma =(\alpha ,\beta )$ where $\alpha\in F_{k-1}$ then  (\ref{2ii}) implies
\begin{eqnarray*}
\int_{(\alpha ,\beta )}~\omega _1~\omega _2~\dots~\omega _{k-1} &=&
\int_{\alpha }~\omega _1~\omega _2~\dots~\omega _{k-1} +\int_{\alpha ^{-1}}~\omega _1~\omega
_2~\dots~\omega _{k-1}\\
& & +\sum_{i=0}^{k-1} \int_{\beta }~\omega _1~\omega
_2~\dots~\omega _{i}\int_{\beta ^{-1} }~\omega _i~\dots~\omega
_{k-1}\\ &=& \int_{\alpha \alpha ^{-1} }~\omega _1~\omega
_2~\dots~\omega _{k-1}+\int_{\beta  ^{-1}\beta  }~\omega _1~\omega
_2~\dots~\omega _{k-1}\\ &=&0 .
\end{eqnarray*}
If, more generally,   $\gamma \in F_{k}$ then  $\gamma =\prod_i \gamma _i$ where each $\gamma $
is   a commutator $(\alpha ,\beta )$, such that
either $\alpha \in F_{k-1}$, or $\beta  \in F_{k-1}$. Therefore
$$
\int_{\gamma }~\omega _1~\omega _2~\dots~\omega _{k-1}=
\sum_i \int_{\gamma _i}~\omega _1~\omega _2~\dots~\omega _{k-1}=0 .
$$
The claim  that $\int_{\gamma }~\omega _1~\omega _2~\dots~\omega
_{k}$, $\gamma \in F_k$ does not depend on the initial point $P_0$
follows from  (iv) (by induction). The claims (iii) and (iv)
follow from (ii) and (\ref{2ii}).

$\Box$

We proved in section \ref{sintegral} that the generating function of limit
cycles $M_k(t)$ is a linear
combination of iterated path integrals of length
$k$ along a loop $\gamma (t)$. As $M(t)$ depends on the {\it free} homotopy class of $\gamma (t)$ then
these iterated integrals are of special nature. The iterated integrals appearing in Lemma \ref{l2} have
the same property : they do not depend on the initial point $P_0$. Therefore they must satisfy (by
analogy to $M_k$) a Fuchsian differential equation. The proof of this  fact can be seen as a
simplified version of the proof of Theorem \ref{main} and for this reason it will be given below.

Let $\gamma (t)\subset f^{-1}(t)$ be a family of closed loops depending continuously on a parameter $t$
in a neighborhood of the typical value $t_0$ of the non-constant polynomial $f\in \C[x,y]$.
We put $S= f^{-1}(t_0)$ and suppose, using the notations of Lemma \ref{l2},  that $\gamma (t_0)\in
F_k$. Consider the
iterated integral
$$ I(t) = \int_{\gamma (t) }~\omega _1~\omega _2~\dots~\omega _{k}
$$ where $\omega _i$ are polynomial one-forms in $\C^2$. In the case $k=1$ this is an Abelian integral
depending on a parameter $t$ and hence it satisfies a (Picard-) Fuchs equation of order
at most $r$.
\begin{proposition} The iterated integral $I$ satisfies a Fuchs equation of order at most $M_r(k)$,
where  $M_r(k) = dim F_k/F_{k+1}$ is given by the Witt formula
\begin{equation}
M_r(k)  = \frac{1}{k} \sum_{d|k} \mu (d) r^{k/d}
\label{mrk}
\end{equation}
and $\mu (d)$ is the M\"obius function (it equals to $0, \pm 1$, see Hall \cite{hall}).
For  small
values of $k,r$ the numbers $M_r(k)$ are shown on the  table below.
\end{proposition}
$$
\begin{tabular}{|l|r|r|r|r|r|r|r|r|}
\hline r $\diagdown$ k &1&2&3&4&5&6&7&8 \\
\hline 1&1&0&0&0&0&0&0&0\\
\hline 2&2&1&2&3&6&9&18&30\\
\hline 3&3&3&8&18&32&116&312&810\\
\hline 4&4&6&20&60&204&4020&4095&8160\\
\hline
\end{tabular}
$$
{\bf Proof.}
Let $\Delta $ be the finite set of atypical values of $f$. The function $I(.)$ is locally analytic on
$\C \setminus \Delta $ and has a moderate growth there (see section \ref{s3}). A finite-dimensional
representation of its monodromy group is constructed as follows.
As $F_k$ is a normal subgroup of $F$ we may consider the
Abelian factor groups
$F_k/F_{k+1}$. Recall that $F_k/F_{k+1}$ is free, torsion free, and finitely
generated.
Thus it is homomorphic to $\Z^{M_r(k)}$ where $r$ is the number of generators of
$F$
and $M_r(k)$ is given by the Witt formula (\ref{mrk}), e.g.
Hall\cite{hall}.
As the Abelian group $F_k/F_{k+1}$ is canonically identified to a subset of $\pi _1(S)$ invariant
under the action (\ref{action1}) of the fundamental group $\pi_1(\C\setminus \Delta  , t_0)$, then we
obtain a homomorphism
\begin{equation}
\pi_1(\C\setminus \Delta  , t_0) \rightarrow Aut(F_k/F_{k+1}) \; .
\label{sub}
\end{equation}
Finally, Lemma 3(ii) implies that the iterated integral $I(t)$ depends on the equivalence class of
$\gamma (t)$ in $F_k/F_{k+1}$. Therefore the monodromy representation of $I$ is a sub-representation of
(\ref{sub}). The Proposition is proved.
\section{ Is the generating function $M_k$ an Abelian integral ? }
Equivalently, is the Fuchs equation satisfied by $M_k$ of Picard-Fuchs type ?
This is an open  difficult problem.
The results of
\cite{gavil} and Theorem \ref{main} provide an answer to the following related question.
Let $f\in \C[x,y]$ be a non constant polynomial, $\gamma (t)\in f^{-1}(t)$ a
family of closed loops
depending continuously on $t$.
Is there a rational one-form on $\C^2$, such that
\begin{equation}
\label{mkt}
M_k(t) = \int_{\gamma (t)} \omega \;\;?
\end{equation}
Consider the  canonical homomorphism
\begin{equation}
\label{hom1} \pi_1 : H_1^\gamma(f^{-1}(t),\Z) \rightarrow
H_1(f^{-1}(t),\Z)
\end{equation}
(which is neither injective, nor surjective in general ) as well
the surjective homomorphism
\begin{equation}\label{hom2}
   \pi_2 : H_1^\gamma(f^{-1}(t),\Z) \rightarrow {\cal M}_k(
\gamma,{\cal F}_\varepsilon ) .
\end{equation}
The homomorphism $\pi_1$ is defined in an obvious way, and $\pi_2$
was defined in section \ref{monrep}. Recall that both of them are
compatible with the action of
 $\pi _1(\C\setminus \Delta , t_0)$.
\\
{\bf Theorem B.} {\it The generating function $M_k$ can be written
in the form (\ref{mkt}) if and only if}
$$
Ker (\pi_1) \subset Ker(\pi_2) .
$$
\\
The theorem says, roughly speaking, that $M_k$ is an Abelian
integral in the sense (\ref{mkt}) if and only if $M_k$ ``depends
on the  homology class of $\gamma (t)$ only". Indeed, when
$M_k(t)$ is an Abelian integral, this holds true. Conversely, if
$Ker (\pi_1) \subset Ker(\pi_2)$, then the injective homomorphism
$$
H_1^\gamma(f^{-1}(t),\Z)/ Ker(\pi_1) \rightarrow H_1(f^{-1}(t),\Z)
$$
and the surjective homomorphism
$$
H_1^\gamma(f^{-1}(t),\Z)/ Ker(\pi_1) \rightarrow {\cal M}_k(
\gamma,{\cal F}_\varepsilon )
$$
are both compatible with the action of $\pi _1(\C\setminus \Delta
, t_0)$. The proof that $M_k$ is an Abelian integral in the sense
(\ref{mkt}) repeats the arguments from the proof of \cite[Theorem
2]{gavil} and will be not
reproduced here.\\
Theorem {\bf B} can be illustrated by the following two basic
examples, taken from \cite{gavil}. \\
 {\bf Example 1.} The
generating function $M_3$ associated to the perturbed foliation
$$
df + \varepsilon(2-x+\frac{1}{2} x^2) dy =0, f= x(y^2-(x-3)^2)
$$
and to the family of ovals $\gamma (t)$ around the center of the
unperturbed system can not be written in the form (\ref{mkt}).
Indeed, an appropriate computation shows that there is a loop
$l(t)$ contained in the orbit of $\gamma (t)$ under the action of
$\pi _1(\C\setminus \Delta , t_0)$, such that
\begin{itemize}
    \item the homology class of $l(t)$ is trivial
    \item the free homotopy class of $l(t)$ is non-trivial
    \item the corresponding generating function $M_3(t)=M_3(l,{\cal
    F}_\varepsilon,t)$ is not identically zero.
\end{itemize}
It follows that $Ker (\pi_1) \not \subset Ker(\pi_2)$ and $M_k(t)$
is not an Abelian integral.\\
 {\bf Example 2.} Let $\omega $ be an
{\it arbitrary} polynomial one-form on $\C^2$. The generating
function $M_k$ associated to
$$
df + \varepsilon \omega  =0, f=y^2 + (x^2-1)^2
$$
and to the exterior family of ovals $\{ f= t\}, t>1$ can be
written in the form (\ref{mkt}). Indeed, it can be shown that the
homomorphism $\pi_1$  (\ref{hom1}) is injective \cite{gavil}. Thus
$M_k$ is always an Abelian integral in the sense (\ref{mkt}),
see also \cite{jmp,jmp1}.

\end{document}